\documentclass[12pt]{amsart}

\usepackage[margin=1in]{geometry}
\usepackage{amsmath,amsthm,caption,dsfont,graphicx,hyperref,makecell,mathtools,multirow,ragged2e,scalerel,tabularx,times,upgreek,amsfonts,xcolor}
\usepackage[capitalise]{cleveref}

\usepackage[foot]{amsaddr}

\newcommand\wh[1]{\widehat{#1}}

\newtheorem{theorem}{Theorem}[section]
\newtheorem{lemma}[theorem]{Lemma}

\newtheorem{remark}[theorem]{Remark}

\newcommand{\R}{\mathbb{R}}
\newcommand{\N}{\mathbb{N}}

\renewcommand{\P}{\mathbf{P}}
\newcommand{\E}{\mathbf{E}}

\newcommand{\eps}{\epsilon}

\DeclareMathSymbol{\shortminus}{\mathbin}{AMSa}{"39}

\title{On the principal eigenvectors of random Markov matrices}

\author{Jacob Calvert}
\address{Institute for Data Engineering and Science, Georgia Institute of Technology, Atlanta, GA 30332}
\email{calvert@gatech.edu}

\author{Frank den Hollander}
\address{Mathematical Institute, Leiden University, Einsteinweg 55, 2333 CC Leiden, The Netherlands}
\email{denholla@math.leidenuniv.nl}

\author{Dana Randall}
\address{School of Computer Science, Georgia Institute of Technology, Atlanta, GA 30332}
\email{randall@cc.gatech.edu}

\begin{document}

\begin{abstract}
We analyze the invariant distributions of continuous-time and discrete-time random walks on randomly weighted complete digraphs. These distributions correspond to the principal left eigenvectors of the associated random Markov generators and kernels, viewed as random matrices. While much is known about the spectra of these matrices, relatively little is known about the principal left eigenvectors, which are delicate random objects for which no explicit form is known. We consider a broad class of such matrices obtained by associating random weights to the vertices and edges of the complete digraph. Our main result concerns the total variation distance between the invariant distribution of the continuous-time random walk and the distribution that is inversely proportional to the vertex weights. It states that, if the edge weights are i.i.d.\ with a finite $p$-th moment for some $p>4$, then this distance a.s.\ converges to zero as the number of vertices grows large, even when the vertex weights are heavy-tailed. We further answer a question of Bordenave, Caputo, and Chafa{\"i} by showing that, despite the dependence of the entries in the corresponding Markov kernel, its invariant distribution is asymptotically uniform a.s., so long as the edge weights have a finite second moment.
\end{abstract}

\maketitle


\section{Introduction and main result}

\noindent
{\bf Background.} Diverse physical and computational processes motivate the study of random Markov matrices, including both the infinitesimal generators of continuous-time Markov chains and the stochastic kernels of discrete-time Markov chains. A natural model of a random Markov generator $Q$ is
\begin{equation}\label{eq:q}
    Q = A - D,
\end{equation}
where $A$ is an $n \times n$ random matrix of non-negative weights that are i.i.d.\ (except possibly for the symmetry $A = A^\top$), and $D = \mathrm{diag}(A 1)$ is the diagonal matrix of the row sums of $A$, ensuring that $Q$ has row sums equal to zero \cite{MR3168123}.\footnote{For a $1 \times n$ vector $v$, we use $\mathrm{diag}(v)$ to denote the $n \times n$ matrix with $v$ along the diagonal and all off-diagonal entries of $0$. For an $n \times n$ matrix $M$, we use $\mathrm{diag}(M)$ to denote the $1 \times n$ vector $(M_{i,i})_{i \in [n]}$.} A closely related model is the random Markov kernel that results from normalizing the rows of $A$ instead:
\begin{equation}\label{eq:k}
    P = D^{-1} A.
\end{equation}
The Markov matrices in \cref{eq:q,eq:k} commonly arise from the adjacency matrices $A$ of graphs, in which case the diagonal entries of $D$ are the vertex degrees and $L = -Q$ is the combinatorial Laplacian \cite{MR3730470}. In this case, $A$ is symmetric when the underlying graph is undirected. In general, $A$ is the weighted adjacency matrix of a directed graph.

Most results about random Markov matrices concern their spectra \cite{MR2575404,MR2644041,MR2857250,MR2892961,MR3168123,MR3189068,MR3717992,MR4440252,MR4716339}. In physical contexts, these spectra provide insight into the relaxation of glassy systems \cite{MR972445}, the resonances of electrical networks \cite{MR1747178,StaringMehligFyodorovLuck2003ImpedanceNetworks}, and other features of the diverse classical and quantum systems that Markov matrices model \cite{Timm2009RandomTransitionRate,MR2539237,MR2488493,Tarnowski2021Superdecoherence,NakerstDenisovHaque2023SparseGenerators}. By comparison, relatively little is known about their principal left eigenvectors, which are invariant distributions $\pi_Q$ and $\pi_P$ of the corresponding continuous-time and discrete-time Markov chains. These eigenvectors relate to the behavior of important algorithms, like PageRank, which ranks the nodes in a graph in terms of the invariant distribution of a random walk running on it \cite{MR3683363,MR4973694}.

By definition, the invariant distributions solve the matrix equations
\begin{equation*}
\pi_Q \, Q = 0 \quad \text{and} \quad \pi_P \, P = \pi_P.
\end{equation*}
When the matrix $A$ in \cref{eq:q,eq:k} is primitive, the Perron--Frobenius theorem implies that $\pi_Q$ and $\pi_P$ exist and are unique. In the special case when $A$ is also Eulerian, or ``balanced'' \cite{MR3730470}, meaning that $A1 = A^\top 1$ (i.e., each row sum equals the corresponding column sum), $\pi_Q$ and $\pi_P$ are simply proportional to the row sums of $A$. This enables the direct characterization of the principal left eigenvectors of symmetric random Markov matrices. For example, if the upper-triangular entries of $A = A^\top$ are i.i.d.\ with finite second moment, then a maximal version of the strong law of large numbers (when applied to the row sums) implies that
\begin{equation}\label{eq: bcc prop}
\| \pi_P - u \|_{\mathrm{TV}} = o(1) \quad \mathrm{a.s.},
\end{equation}
where $u = (\frac1n, \dots, \frac1n)$ is the uniform distribution and $\|\nu\|_{\mathrm{TV}} = \frac12 \sum_{i \in [n]} |\nu (i)|$ denotes the total variation norm \cite[Proposition 1.5]{MR2644041}.

In general, however, the principal left eigenvector is a ``delicate random object with no explicit expression'' \cite{MR3916108}. Indeed, with $M$ denoting $Q$ or $P$, the Markov chain tree theorem states that
\begin{equation*}
\pi_M (i) \propto \sum_{\mathcal{T}_i} \prod_{(j,k) \in \mathcal{T}_i} M_{j,k},
\end{equation*}
where the sum ranges over spanning trees $\mathcal{T}_i$ of the complete digraph rooted at vertex $i$ \cite{MR734725}. The complexity of this expression gives a sense of the challenge that the characterization of $\pi_M$ generally poses.

\medskip\noindent
{\bf Goal.} The main purpose of this paper is to show that, although $\pi_Q$ is generally a complicated function of the transition rates, for a broad class of random Markov generators it is closely approximated by a simple function of the exit rates
\begin{equation*}
q_i = |Q_{i,i}| = \sum_{j:\, j \neq i} Q_{i,j}.
\end{equation*}
Regarding $\pi_P$, we extend \cref{eq: bcc prop} to the general case of non-symmetric adjacency matrices $A$, thereby addressing several related questions and claims in the literature.

\medskip\noindent
{\bf Model.} We study the Markov matrices $Q$ and $P$ in \cref{eq:q,eq:k} arising from an $n \times n$ random weighted adjacency matrix~$A$ with entries
\begin{equation*}
A_{i,j} = \theta_i X_{i,j}, \quad i, j \in [n],
\end{equation*}
where $(\theta_i)_{i \geq 1} \in (0,\infty)^\N$ is a (possibly random) sequence of vertex weights and $(X_{i,j})_{i,j \geq 1}$ is an infinite array of i.i.d.\ directed edge weights with non-degenerate law $\mathcal{L}$ supported on $[0,\infty)$.\footnote{This (non-unique) parameterization of adjacency matrices is common in stochastic thermodynamics and chemical kinetics, where the logarithms of $\theta_i$ and $X_{i,j}$ represent the well depths and effective barrier heights in the energy landscape of a non-equilibrium system \cite{PhysRevX.10.011066}.} Note that the vertex weights do not affect $P$, nor the closely related kernel $\wh{Q}$ defined by
\[
\wh{Q} = \wh{D}^{-1} \wh{A}
\]
in terms of $\wh{A} = A - \mathrm{diag}(\mathrm{diag} (A))$ and $\wh{D} = \mathrm{diag}(\wh{A} 1)$. We will not study $\wh{Q}$ for its own sake. Rather, its analysis will be key to proving our main result concerning $\pi_Q$.

Although $\mathcal{L}$ may have an atom at $0$, the edge weight $X_{i,j}$ is positive with some probability $p \in (0,1)$. It follows that both $A$ and $\wh{A}$ are primitive for all sufficiently large~$n$ a.s., because the same is true of the adjacency matrix of the directed Erd\H{o}s--R\'{e}nyi random graph $\mathcal{D}(n,p)$ with and without self-loops. For ease of reference, we collect some standard yet important consequences of this observation in the following remark.

\begin{remark}\label{rem: A prim}
    For the rest of the paper, we will work on the event that $A$ and $\wh{A}$ are primitive, which holds for all sufficiently large $n$ a.s. On this event, the Markov matrices $Q$, $P$, and $\wh{Q}$ are well-defined and irreducible, and hence the corresponding invariant distributions $\pi_Q$, $\pi_P$, and $\pi_{\wh Q}$ exist and are unique. Additionally, the kernels $P$ and $\wh{Q}$ are aperiodic.
\end{remark}

\begin{figure}
    \centering
    \includegraphics[width=\linewidth,trim={0 1cm 0 0.5cm},clip]{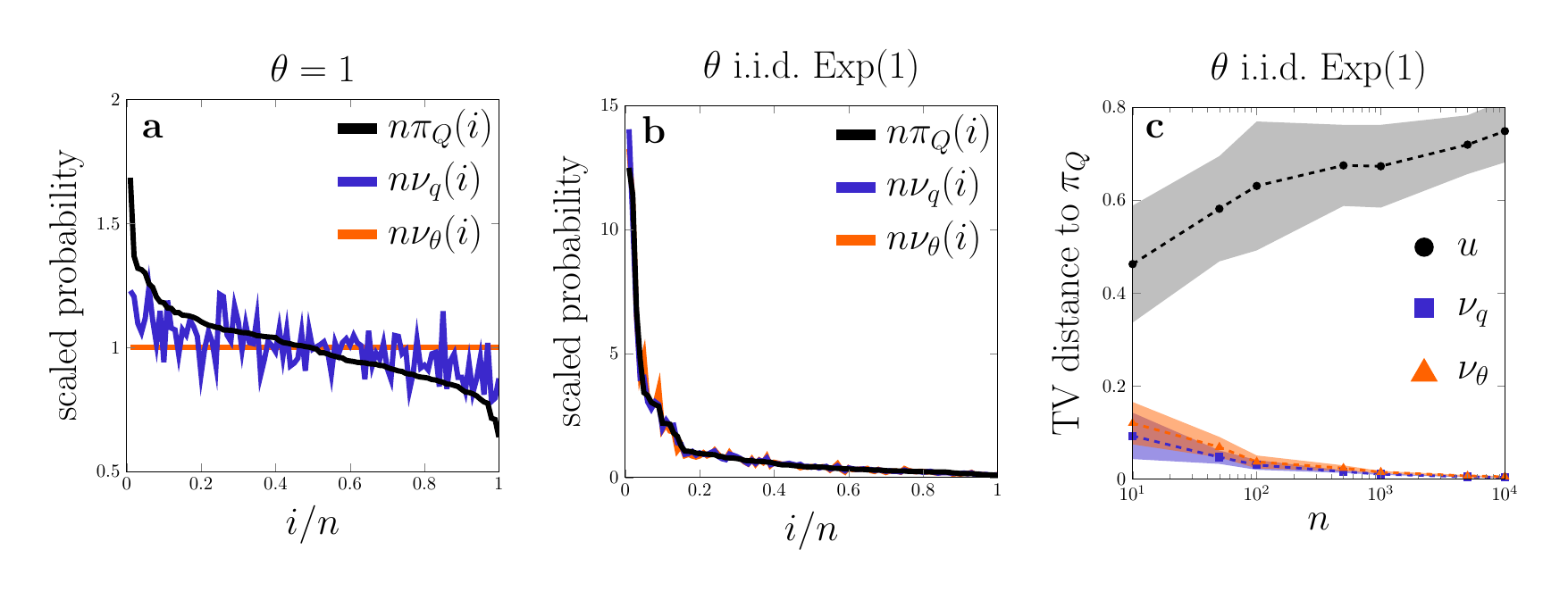}
    \caption{Comparison of distributions in \cref{thm:piq} for edge weight law $\mathcal{L} = \mathrm{Exp}(1)$. Plotted are the scaled $\pi_Q$ probabilities in descending order with corresponding values of $\nu_q$, $\nu_\theta$ for $n=100$ and vertex weights $\theta_i$: (\textbf{a}) all equal to $1$; (\textbf{b}) independent $\mathrm{Exp}(1)$. In (\textbf{c}), as $n$ increases, $\pi_Q$ grows less uniform and closer to $\nu_q$, $\nu_\theta$ in total variation (TV) distance. In all subplots, lines represent averages over $10$ trials and shaded regions represent $\pm 1$ standard deviation.
    }
    \label{fig: summary}
\end{figure}

\medskip\noindent
{\bf Main results.} Our first result states that $\pi_Q$ is asymptotically inversely proportional to the exit rates and vertex weights whenever the edge weights have a finite $p$-th moment for some $p > 4$. For $x \in (0,\infty)^n$, we denote by $\nu_x$ the distribution that satisfies $\nu_x (i) \propto x_i^{-1}$.

\begin{theorem}\label{thm:piq}
	If $\mathcal{L}$ has a finite $p$-th moment for some $p > 4$, then, for every $a < \min\{\frac12,1 - \frac{4}{p}\}$,
	\begin{equation*}
		\left\| \pi_Q - \nu \right\|_{\mathrm{TV}} = O (n^{-a}) \quad \mathrm{a.s.}, \quad\quad \nu \in \{\nu_q, \, \nu_\theta\}.
	\end{equation*}
\end{theorem}

Because \cref{thm:piq} places no restriction on the vertex weights, which can be deterministic (\cref{fig: summary}a) or random (\cref{fig: summary}b), the distributions $\pi_Q$, $\nu_q$, and $\nu_\theta$ can all be far from uniform. For example, when all vertex and edge weights are i.i.d.\ $\mathrm{Exp}(1)$ random variables, these distributions grow closer in total variation distance---and farther from uniform---as $n$ increases (\cref{fig: summary}c). Moreover, when the entries of $\theta$ are i.i.d.\ with a heavy-tailed distribution, \cref{thm:piq} implies that $\pi_Q$ localizes on ``trap'' states with extremely low exit rates or vertex weights \cite{MR2581889}. For comparison, in the reversible case when $\theta \equiv 1$ and the edges are undirected, $\pi_P (i)$ is known to be exactly proportional to $\sum_j X_{i,j}$, in which case i.i.d.\ heavy-tailed edge weights instead produce traps along the edges of highest weight \cite[Theorem 1.8]{MR2857250}.

The intuition behind \cref{thm:piq} is that, although the generator $Q$ is random, the transition probabilities $\wh{Q}_{i,j} = Q_{i,j}/q_i$ with which the embedded discrete-time Markov chain jumps to neighboring states are comparable. As a consequence, the latter Markov chain rapidly converges to a stationary distribution $\pi_{\wh{Q}}$ that is close to uniform. Since $\pi_Q (i)$ is proportional to $\pi_{\widehat{Q}} (i)/q_i$, $\pi_Q$ is essentially determined by the reciprocal exit rates. The exit rates $q_i$, in turn, are approximately $\theta_i n\, \E (X_{1,1})$, due to the concentration of sums of edge weights.

\begin{remark}\label{rem: tail}
In \cref{thm:piq}, we opt for a simple moment assumption on $\mathcal{L}$, but it is possible to prove qualitatively similar results under a technically sharper tail condition of the form
\[
\sum_{k \geq 1} n_k^2 \, \bar{F} \big(\sqrt{n_k}/g(n_k)\big) < \infty,
\]
where $\bar{F}(x) = \mathcal{L}([x,\infty))$ is the complementary cumulative distribution function, $g: \N \to (0,\infty)$ is a function satisfying $g(n) = \Omega (\sqrt{\log n}\,)$, and $(n_k)_{k \geq 1} \subseteq \N$ is an increasing sequence. For example, taking $n_k = 2^k$ and $g(n) = (\log n) \vee 1$ results in a condition that is stronger than $\mathcal{L}$ having a finite fourth moment, but weaker than the finiteness of a $p$-th moment for some $p > 4$. In this case, the total variation distance in \cref{thm:piq} is $O ((\log n)^{-1})$ instead.
\end{remark}

The key input to the proof of \cref{thm:piq} is an $o(1/n)$ estimate of the $\ell^\infty$ distance of $\pi_{\widehat{Q}}$ to the uniform distribution $u$, which uses the finiteness of a $p$-th moment for some $p > 4$. However, our second result states that it suffices for the edge weights to have a finite second moment to guarantee the asymptotic uniformity of $\pi_{\wh{Q}}$ in the $\ell^1$ sense. The same is true of $\pi_P$.

\begin{theorem}\label{thm:m2}
	If $\mathcal{L}$ has a finite second moment, then 
	\begin{equation}\label{eq:m2 thm}
		\left\| \pi - u \right\|_{\mathrm{TV}} = o(1) \quad \mathrm{a.s.}, \quad\quad \pi \in \{\pi_P,\pi_{\wh{Q}}\}. 
	\end{equation}
	Furthermore, if $\theta \equiv c \in (0,\infty)$, then \cref{eq:m2 thm} also holds with $\pi = \pi_Q$.
\end{theorem}

\begin{figure}
    \centering
    \includegraphics[width=\linewidth,trim={0 1cm 0 0.5cm},clip]{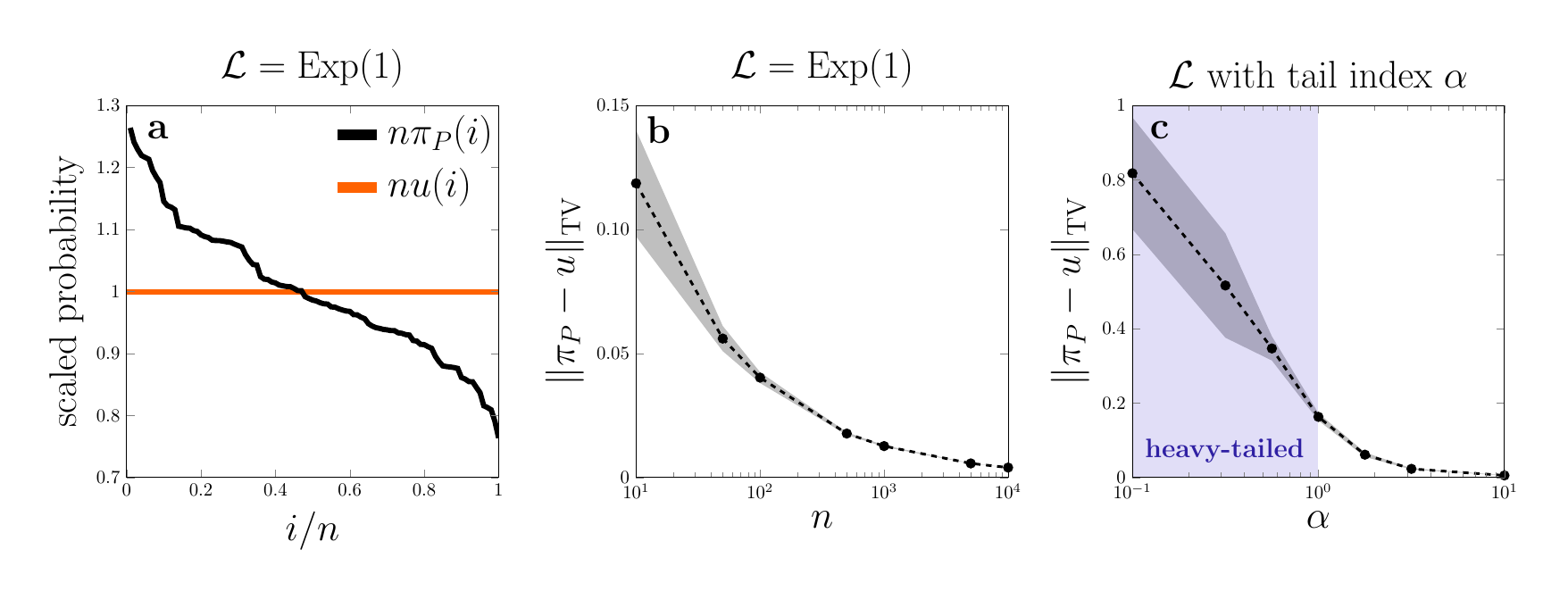}
    \caption{Uniformity of $\pi_P$. (\textbf{a}) Plotted are the scaled $\pi_P$ probabilities in descending order for $n=100$. Plotted are also the TV distance of $\pi_P$ to uniform with: (\textbf{b}) $\mathrm{Exp}(1)$ edge weights, as $n$ varies (note $y$-axis upper limit); (\textbf{c}) independent $X_{i,j} \stackrel{\mathrm{d}}{=} Y^{-1/\alpha}$ edge weights, where $Y \sim \mathrm{Exp}(1)$, as $\alpha > 0$ varies for $n = 100$ fixed. In (\textbf{c}), $\mathcal{L}$ has infinite mean when $\alpha \leq 1$.
    Lines and shading reflect averages and $\pm 1$ standard deviation over $10$ trials.
    }
    \label{fig: summary alphas}
\end{figure}

The conclusion of \cref{thm:m2} concerning $\pi_P$ answers a question of Bordenave, Caputo, and Chafa{\"i} \cite{MR2892961}. This follows Chafa{\"i}'s earlier suggestion to study the asymptotic behavior of $\pi_P$ for $P$ from the Dirichlet Markov ensemble \cite{MR2575404}, meaning that the rows of $P$ are independent $\mathrm{Dirichlet}(u)$ random variables, or equivalently that $\mathcal{L} = \mathrm{Exp}(1)$ (\cref{fig: summary alphas}a). This ensemble of random Markov matrices is particularly common in applications, as it is the uniform or maximal entropy distribution over such matrices. \cref{thm:m2} implies the asymptotic uniformity of such $\pi_P$ (\cref{fig: summary alphas}b), and verifies claims of Horvat \cite{MR2539237} and Edelman et al.\ \cite{edelman2024the}.

We further note that \cite{MR3168123} studied the same model but without vertex weights, i.e., $\theta = 1$ (\cref{fig: summary}a). Under the assumption that $\mathcal{L}$ has a finite fourth moment, they proved that
\[
\| \pi_Q - u \|_{\mathrm{TV}} = O \big( \sqrt{(\log n)/n} \,\big) \quad \mathrm{a.s.}
\]
\cref{thm:m2} establishes that $\pi_Q$ is asymptotically uniform a.s.\ under the weaker assumption of edge weights with a finite second moment. 

\medskip\noindent
{\bf Open problems.} What conditions on $\mathcal{L}$ are necessary for results like \cref{thm:piq,thm:m2} to hold? To significantly weaken the condition in \cref{thm:piq} beyond what \cref{rem: tail} discusses would require a new proof idea. However, even when $\mathcal{L}$ is heavy-tailed, $\nu_q$ may closely approximate $\pi_Q$ in the sense of $\log \pi_Q (i)$ and $\log \nu_q (i)$ having correlation close to $1$ for a random state $i$ \cite{MR4826995,MR5019616}. \cref{fig: summary alphas}c shows that $\pi_P$ can be far from uniform when $\mathcal{L}$ is heavy-tailed. Is it necessary for $\mathcal{L}$ to have a finite second moment for \cref{eq:m2 thm} to hold with $\pi = \pi_P$, as it is when $P$ is symmetric \cite{MR2857250}? 

\medskip\noindent
{\bf Organization.} In \cref{sec: estimates}, we prove bounds on transition probabilities that serve as inputs to the later sections. \cref{sec: proof of thm1} establishes a key estimate of the $\ell^\infty$ uniformity of $\pi_{\wh{Q}}$, which we then use to prove \cref{thm:piq}. Lastly, in \cref{sec: proof of thm2}, we prove \cref{thm:m2} by obtaining lower bounds on the entries of $P^2$ and $\wh{Q}^2$. 


\section{Estimates of transition probabilities}\label{sec: estimates}

We will repeatedly use a maximal version of the strong law of large numbers due to Bai and Yin \cite[Lemma 2]{MR1235416}.

\begin{lemma}\label{lem:bai yin}
Let $(X_{i,j})_{i, j \geq 1}$ be an array of i.i.d.\ random variables. For any real numbers $\alpha \in (1/2,1]$, $\beta \geq 0$, and $C > 0$, if $\E(|X_{1,1}|^{(1+\beta)/\alpha}) < \infty$, then
\[
		\max_{i \in [Cn^\beta]} \left| \sum_{j \in [n]} (X_{i,j} - \E (X_{1,1})) \right| = o(n^{\alpha}) \quad \mathrm{a.s.}
\]
\end{lemma}

For the proof of \cref{thm:m2}, we will also need to control the minimum among $n^2$ sums of elements from $X$. Under the assumption of finite second moment alone, \cref{lem:bai yin} is not enough. However, because the entries of $X$ are non-negative, the lower tail of their sums is exponential.

\begin{lemma}\label{lem:lower tail chernoff}
Let $(Y_i)_{i \geq 1}$ be a sequence of i.i.d.\ non-negative random variables with mean $\mu$ and variance $\sigma^2 \in (0,\infty)$. For every integer $n \geq 1$ and $\eps \in (0,1)$, the sum $S_n = \sum_{i \in [n]} Y_i$ satisfies
\begin{equation}\label{eq: sn chernoff}
    \P (S_n \leq (1-\eps) \mu n) \leq \exp \left(- \frac{\eps^2 \mu^2 n}{2(\mu^2 + \sigma^2) } \right).
\end{equation}
\end{lemma}

\begin{proof}
For every $t > 0$, the application of Markov's inequality to $e^{-t S_n} \geq 0$ implies that
\begin{equation}\label{eq: sn markov}
	\P (S_n \leq (1-\eps)\mu n) = \P \left(e^{-t S_n} \geq e^{-t (1-\eps) \mu n} \right) \leq e^{t(1-\eps) \mu n}\, \E (e^{-t S_n}).
\end{equation}
By the independence of the $Y_i$ and the basic inequalities $e^{-x} \leq 1 - x + \frac12 x^2$ and $1+y \leq e^y$, which hold for all $x \geq 0$ and $y \in \R$, the expectation satisfies
\[
	\E (e^{-t S_n}) = \big( \E (e^{-tY_1}) \big)^n \leq \big(1 - t \mu + \tfrac12 t^2 m_2 \big)^n \leq \exp \left(- (t \mu - \tfrac12 t^2 m_2)n \right)	
\]
with $m_2 = \E(Y_1^2) = \mu^2 + \sigma^2$. Substitution into \cref{eq: sn markov} yields the bound
\[
	\P (S_n \leq (1-\eps) \mu n) \leq \exp \left( - (t \eps \mu - \tfrac12 t^2 m_2 )n \right),
\]
which is minimized at $t = \eps \mu /m_2$, and gives \cref{eq: sn chernoff}.
\end{proof}

\cref{lem:bai yin} implies the following estimate of sums of squared entries of the kernels $\wh{Q}$ and $P$.

\begin{lemma}\label{lem:l2 control}
Let $K$ denote either $\wh{Q}$ or $P$. If $\mathcal{L}$ has finite fourth moment, then
\[
	\max_{i \in [n]} \sum_{j \in [n]} K_{i,j}^{\,2} = O(n^{-1}) \quad \mathrm{a.s.}
\]
\end{lemma}

\begin{proof}
Denote the $i$-th off-diagonal row sum of $X$ by $R_i = \sum_{j: \, j \neq i} X_{i,j}$. By \cref{lem:bai yin}, because $X_{1,1}$ has finite second moment, these row sums satisfy
\[
	\min_{i \in [n]} R_i = (1- o(1))\, \mu n \quad \mathrm{a.s.}
\]
with $\mu = \E(X_{1,1})$. We further apply \cref{lem:bai yin} to the array $(Y_{i,j})_{i,j \geq 1}$ of squared entries $Y_{i,j} = X_{i,j}^2$. Because $\E (Y_{1,1}^2) = \E (X_{1,1}^4) < \infty$, \cref{lem:bai yin} implies that
\[
	\max_{i \in [n]} \sum_{j: \, j \neq i} Y_{i,j} = (1+o(1))\, \E (X_{1,1}^2)\, n \quad \mathrm{a.s.}
\]
On these two a.s.\ events, we have that
\[
	\max_{i \in [n]} \sum_{j \in [n]} K_{i,j}^{\, 2} \leq \frac{\max_{i \in [n]} \sum_{j \in [n]} Y_{i,j}}{(\min_{i \in [n]} R_i)^2} \leq \frac{(1+o(1))\, \E (X_{1,1}^2)\, n}{( (1-o(1)) \mu n )^2} = O(n^{-1}).
\]
Note that the first inequality holds because the numerator of the upper bound includes diagonal entries of $X$ that $\wh{Q}$ excludes, while the row sum $R_i$ in the denominator excludes the diagonal entries of $X$ that $K = P$ includes.
\end{proof}

We use the finiteness of a $p$-th moment for some $p > 4$ to control the largest entries of the kernels.

\begin{lemma}\label{lem:max control}
Let $K$ denote $\wh{Q}$ or $P$. If $\mathcal{L}$ has finite $p$-th moment for some $p >4$, then, for every $a < 1 - \frac{2}{p}$,
\[
	\max_{i,j \in [n]} K_{i,j} = O (n^{-a}) \quad \mathrm{a.s.}
\]
\end{lemma}

\begin{proof}
By \cref{lem:bai yin}, since $X_{1,1}$ has finite second moment, the largest row sums of $\wh{Q}$ and $P$ are both $\Omega (n)$ a.s. On this event, it suffices to show that
\[
	W_n = \max_{i,j \in [n]} X_{i,j} = O(n^{1-a}) \quad \mathrm{a.s.}
\]
In fact, since $W_n$ is increasing in $n$, it is enough for this event to hold along the dyadic subsequence $n_k = 2^k$. Let $\delta = 1 - \frac{2}{p} - a > 0$. A union bound over the $n_k^2$ possible values of $i$ and $j$ and the application of Markov's inequality to $X_{1,1}^p$ together imply that
\[
    \P (W_{n_k} > n_k^{1-a}) = O \big(n_k^2 n_k^{-p(1-a)} \big) = O \big(n_k^{-p\delta} \big) = O\big(2^{-p\delta k}\big).
\]
Because this bound is summable over $k \geq 1$, the Borel--Cantelli lemma implies that $W_{n_k} = O(n_k^{1-a})$ a.s.
\end{proof}


\section{Proof of \texorpdfstring{\cref{thm:piq}}{Theorem 1.2}}\label{sec: proof of thm1}

We next address the key input to the proof of \cref{thm:piq}: an $o(n^{-1})$ estimate of the $\ell^\infty$ distance between $\pi_{\wh{Q}}$ and the uniform distribution $u$.

\begin{lemma}\label{lemma:P l_infty}
If $\mathcal{L}$ has finite $p$-th moment for some $p > 4$, then, for every $a < \min\{\frac12,1-\frac{4}{p}\}$,
\begin{equation}\label{eq:P l_infty}
	\| \pi_{\wh{Q}} - u \|_\infty = O(n^{-(1+a)}) \quad \mathrm{a.s.}
\end{equation}
\end{lemma}

The idea of the proof is to use Bernstein's inequality to obtain exponential concentration of the two-step transition kernel $(\wh{Q}^2)_{i,k}$ given the $i$-th row of $\wh{Q}$. The exponential rate of concentration is determined by a balance of the $\ell^2$ and $\ell^\infty$ norms of the $i$-th row of $\wh{Q}^2$, for which \cref{lem:l2 control,lem:max control} provide a.s.\ estimates.

\begin{proof}[Proof of \cref{lemma:P l_infty}]
	Fix $p > 4$ and $a < \min\{\frac12,1-\frac{4}{p}\}$. According to \cref{rem: A prim}, $\wh{Q}$ is irreducible and aperiodic, and hence $\pi_{\wh{Q}} = \pi_{\wh{Q}^2}$ for all sufficiently large $n$ a.s. It therefore suffices to show that
\begin{equation}\label{eq:p2ik lb}
	\max_{i,k \in [n]} \left| (\wh{Q}^2)_{i,k} - \frac1n \right| = O\big(n^{-(1+a)} \big) \quad \mathrm{a.s.}
\end{equation}
Indeed, because $\pi_{\wh{Q}}$ is stationary, on the event in \cref{eq:p2ik lb}, it satisfies
\[
	\pi_{\wh{Q}} (i) = \sum_{j \in [n]} \pi_{\wh{Q}} (j) (\wh{Q}^2)_{j,i} = \frac{1}{n} + O(n^{-(1+a)}),
\]
which implies \cref{eq:P l_infty}.

To prove \cref{eq:p2ik lb}, we begin by noting that, conditionally on the $i$-th row $\wh{Q}_{i,\cdot}$ of $\wh{Q}$, the random variables $Y_j = \wh{Q}_{i,j} (\wh{Q}_{j,k} - \frac{1}{n-1})$ indexed by $j \in [n]\setminus\{i,k\}$ are independent, centered, and satisfy
\begin{equation}\label{eq:sumy}
	S_n (i,k) = \sum_{j \in [n]\setminus\{i,k\}} Y_j = (\wh{Q}^2)_{i,k} - \frac{1-\wh{Q}_{i,k}}{n-1}. 
\end{equation}
Bernstein's inequality states that there is a universal constant $C > 0$ such that, for every $t > 0$,
\begin{equation}\label{eq:bernstein}
	\P \left( |S_n (i,k)| > t \mid \wh{Q}_{i,\cdot} \right) \leq 2 \exp \left(- C \min\left\{\frac{t^2}{v_n}, \frac{t}{b_n} \right\} \right),
\end{equation}
with
\[
	v_n = \sum_{j \in [n]\setminus\{i,k\}} \E (Y_j^2 \mid \wh{Q}_{i,\cdot}), \qquad b_n = \max_{j \in [n]\setminus\{i,k\}} |Y_j|.
\]
Regarding $v_n$, \cref{lem:bai yin} implies that, because $\mathcal{L}$ has finite second moment, the off-diagonal row sums $R_i = \sum_{j:\, j \neq i} X_{i,j}$ satisfy
\[
	\min_{i \in [n]} R_i = (1-o(1))\, \E(X_{1,1}) n \quad \mathrm{a.s.}
\]
On this event, as well as the a.s.\ event of \cref{lem:l2 control}, the fact that $X_{j,k}$ is independent of $\wh{Q}_{i,\cdot}$ and has finite second moment implies that
\[
	v_n = \sum_{j \in [n]\setminus\{i,k\}} \wh{Q}^2_{i,j} \, \E (\wh{Q}^2_{j,k} \mid \wh{Q}_{i,\cdot}) = O(n^{-2}) \sum_{j \in [n]\setminus\{i,k\}} \wh{Q}^2_{i,j} = O(n^{-3}).
\]
Furthermore, on the a.s.\ event of \cref{lem:max control}, for every $s < 2(1 - \frac{2}{p})$, 
\[
	b_n \leq \max_{j \in [n]\setminus\{i,k\}} \wh{Q}_{i,j} \Big(\wh{Q}_{j,k} + \frac{1}{n-1}\Big) = O(n^{-s}).
\]
In particular, we choose $s > 1+a$.

When $t_n = \Theta (n^{-(1+a)})$ the rate in \cref{eq:bernstein} satisfies
\[
	\min\left\{\frac{t_n^2}{v_n}, \frac{t_n}{b_n} \right\} = \Omega \left(\min \left\{ n^{1-2a}, n^{s-1-a}\right\}\right) = \Omega (n^{\delta}) \quad \mathrm{a.s.},
\]
where $\delta = \min\{1-2a,s-1-a\}> 0$. A union bound over the $n^2$ possible values of $i$ and $k$, together with \cref{eq:bernstein}, imply that
\[
	\P \left( \max_{i,k \in [n]} |S_n (i,k)| > t_n \right) \leq 2n^2 \exp \big(- \Omega(n^\delta) \big) \quad \mathrm{a.s.}
\]
Because this upper bound is summable, the Borel--Cantelli lemma implies that
\[
	\max_{i,k \in [n]} |S_n (i,k)| = O (n^{-(1+a)}) \quad \mathrm{a.s.}
\]

To relate this to \cref{eq:p2ik lb}, we combine the preceding fact with \cref{eq:sumy}. The triangle inequality and \cref{lem:max control} imply that, for every $c < 2(1-\frac{1}{p})$,
\[
	\max_{i,k \in [n]} \left| (\wh{Q}^2)_{i,k} - \frac1n \right| \leq \max_{i,k \in [n]} \left( \left| S_n (i,k) \right| + \left| \frac{1-\wh{Q}_{i,k}}{n-1} - \frac1n \right| \right) = O (n^{-(1+a)}) + O(n^{-c}) \quad \mathrm{a.s.},
\]
In particular, we can take $c = \frac{3}{2} > 1 + a$ because $p > 4$, in which case $O(n^{-c}) = o(n^{-(1+a)})$. This proves \cref{eq:p2ik lb}, and completes the proof.
\end{proof}

\cref{thm:piq} follows readily from the fact that $\pi_Q(i) \propto \pi_{\wh{Q}} (i)/q_i$ and \cref{lemma:P l_infty}.

\begin{proof}[Proof of \cref{thm:piq}]
	Put $\delta_i = n \pi_{\wh{Q}} (i) - 1$. Since $\pi_Q (i)$ is proportional to $\pi_{\wh{Q}} (i) / q_i$, it satisfies
	\[
		\pi_Q (i) = \frac{\pi_{\wh{Q}} (i) / q_i}{\sum_{k \in [n]} \pi_{\wh{Q}} (k) / q_k} = \frac{(1 + \delta_i) / q_i}{\sum_{k \in [n]} (1 + \delta_k) / q_k}.
	\]
	In terms of $\delta = \max_{i \in [n]} |\delta_i|$, this implies that
	\[
		\left( \frac{1-\delta}{1+\delta} \right) \nu_q (i) \leq \pi_Q (i) \leq \left( \frac{1+\delta}{1-\delta}\right) \nu_q (i).
	\]
	Since $\delta = n\|\pi_{\wh{Q}} - u\|_\infty$, \cref{lemma:P l_infty} implies that, for every $a < \min\{\frac12,1-\frac{4}{p}\}$, $\delta = O(n^{-a})$ a.s. On this event, the preceding bound implies that
	\[
		\| \pi_Q - \nu_q \|_{\mathrm{TV}} \leq \max \left\{\frac{\delta}{1 + \delta}, \frac{\delta}{|1 - \delta|} \right\} = O (n^{-a}).
	\]
    The fact that we obtain the same estimate with $\nu_\theta$ in place of $\nu_q$ is a consequence of the triangle inequality and the following lemma.
\end{proof}

\begin{lemma}
	If $\mathcal{L}$ has a finite $p$-th moment for some $p > 4$, then, for every $a < \min\{\tfrac12,1-\tfrac{3}{p}\}$,
	\begin{equation*}
	\| \nu_q - \nu_\theta \|_{\mathrm{TV}} = O(n^{-a}) \quad \mathrm{a.s.}
	\end{equation*}
\end{lemma}

\begin{proof}
Let $\mu = \E(X_{1,1})$. Because the $Y_{i,j} = X_{i,j} - \mu$ are independent mean-zero random variables with finite $p$-th absolute moment for $p > 2$, the Fuk--Nagaev inequality (see, e.g., \cite[Eq.\ 1.8]{MR3652041}) states that there are constants $C_1$ and $C_2$ depending on $\mathcal{L}$ such that, for any $z > 1$,
\begin{equation*}
	\P \left( \Big| \sum_{j:\, j \neq i} Y_{i,j} \Big| \geq C_1 \sqrt{n \log z} + C_2 (n z)^{1/p} \right) \leq 4/z.
\end{equation*}
Take $z = (n \log n)^2$. Combining a union bound over $i \in [n]$ with the preceding inequality, yields, for all sufficiently large $n$,
\[
	\P \left( \max_{i \in [n]} \,\Big| \sum_{j:\, j \neq i} Y_{i,j} \Big| \geq C_1' \sqrt{n \log n} + C_2 n^{3/p} (\log n)^{2/p}\right) \leq \frac{4}{n (\log n)^2}
\]
with $C_1'$ a further constant. Because this upper bound is summable, the Borel--Cantelli lemma implies that 
\[
	\max_{i \in [n]} \Big| \sum_{j:\, j \neq i} X_{i,j} - (n-1) \mu \Big| = \max_{i \in [n]} \,\Big| \sum_{j:\, j \neq i} Y_{i,j} \Big| = O \left( \sqrt{n \log n} + n^{3/p}(\log n)^{2/p} \right) \quad \mathrm{a.s.}
\]
In particular, for every $a < \min\{\tfrac12, 1 - \tfrac{3}{p}\}$,
\[
	\max_{i \in [n]} \frac{1}{n\mu} \Big| \sum_{j:\, j \neq i} X_{i,j} - n \mu \Big| = O(n^{-a}) \quad \mathrm{a.s.}
\]
Recall that $q_i = \theta_i \sum_{j:\, j \neq i} X_{i,j}$. Hence, on the preceding event, $\nu_q$ satisfies
\[
	\nu_q (i) = \frac{1/q_i}{\sum_{k \in [n]} 1/q_k} = \frac{1/(\theta_i n \mu)}{\sum_{k \in [n]} 1/(\theta_k n \mu)} \left(1 + O(n^{-a}) \right) = \nu_\theta (i) \left(1 + O(n^{-a}) \right),
\]
which implies that
\[
	\| \nu_q - \nu_\theta \|_{\mathrm{TV}} = \frac12 \sum_{i \in [n]} \nu_\theta (i) \, O(n^{-a}) = O(n^{-a}).
\]
\end{proof}


\section{Proof of \texorpdfstring{\cref{thm:m2}}{Theorem 1.4}}\label{sec: proof of thm2}

As in the proof of \cref{lemma:P l_infty}, we establish \cref{thm:m2} by considering the two-step transition matrices. However, under the assumption of finite second moment alone, we cannot obtain a two-sided estimate like \cref{eq:p2ik lb} for their entries. Instead, we merely bound from below their smallest entries by using the exponential concentration of the lower tail of sums of $X$ (\cref{lem:lower tail chernoff}). This is enough to bound the $\ell^1$ norms in \cref{eq:m2 thm}, because their row sums are normalized.

\begin{proof}[Proof of \cref{thm:m2}]

Let $K$ denote either $P$ or $\wh{Q}$. Following \cref{rem: A prim}, $\pi_K = \pi_{K^2}$ for all sufficiently large $n$ a.s. It therefore suffices to show that
\begin{equation}\label{eq:mik lb}
	\min_{i,k \in [n]} (K^2)_{i,k} \geq \frac{1 - o(1)}{n} \quad \mathrm{a.s.}
\end{equation}
Indeed, because $\pi_{M}$ is stationary, \cref{eq:mik lb} implies that
\[
	\pi_{K} (i) = \sum_{j \in [n]} \pi_{K} (j) (K^2)_{i,j} \geq \frac{1 - o(1)}{n} \quad \mathrm{a.s.}
\]
On this event, defining $\mathcal{A} = \{i \in [n]: \pi_{K} (i) \leq 1/n\}$, we have that
\[
\left\| \pi_{K} - u \right\|_{\mathrm{TV}} = \frac12 \left(\frac{|\mathcal{A}|}{n} - \pi_{K} (\mathcal{A}) + \pi_{K} (\mathcal{A}^c) - \frac{|\mathcal{A}^c|}{n}\right) = \frac{|\mathcal{A}|}{n} - \pi_{K} (\mathcal{A}) = o \left( \frac{|\mathcal{A}|}{n} \right) = o(1).
\]

To prove \cref{eq:mik lb}, note that
\begin{equation}\label{eq:mik interm lower bound}
	(K^2)_{i,k} \geq \frac{S_n (i,k)}{(\max_{i \in [n]} R_i)^2},
\end{equation}
where $R_i = \sum_{j \in [n]} X_{i,j}$ and $S_n (i,k) = \sum_{j \in [n]\setminus\{i,k\}} X_{i,j} X_{j,k}$. By \cref{lem:bai yin}, since $X_{1,1}$ has finite second moment, the row sums satisfy
\begin{equation}\label{eq:col sum ub}
	\max_{i \in [n]} R_i = (1 + o(1))\, \mu n \quad \mathrm{a.s.}
\end{equation} 
with $\mu = \E(X_{1,1})$. Additionally, because the summands of $S_n (i,k)$ are i.i.d.\ random variables with finite variance, \cref{lem:lower tail chernoff} implies that there is a constant $c > 0$ such that, if $\eps_n = c( (\log n)/n )^{1/2}$, then
\[
	\P ( S_n (i,k) \geq (1-\eps_n)\, \mu^2 n ) \leq n^{-4}.
\]
As the union bound over the $n^2$ possible values of $i$ and $k$ is summable, the Borel--Cantelli lemma implies that 
\begin{equation}\label{eq:numerator lower bd}
	\min_{i,k \in [n]} S_n (i,k) \geq (1 - o(1))\, \mu^2 n \quad \mathrm{a.s.}
\end{equation}
By \cref{eq:mik interm lower bound}, on the combined events of \cref{eq:col sum ub,eq:numerator lower bd}, $K^2$ satisfies
\[
	\min_{i,k \in [n]} (K^2)_{i,k} \geq \frac{(1-o(1))\, \mu^2 n}{((1+o(1))\, \mu n)^2} = \frac{1-o(1)}{n},
\]
which establishes \cref{eq:mik lb}, and completes the proof of \cref{eq:m2 thm} for $\pi_P$ and $\pi_{\wh{Q}}$.

Regarding $\pi_Q$, when $\theta \equiv c \in (0,\infty)$, the row sums of $Q$ are equal to those of $X$. Because $\mathcal{L}$ has finite second moment, \cref{lem:bai yin} implies that the largest $\delta_i = \mu n / q_i - 1$ satisfies
\[
	\max_{i \in [n]} |\delta_i| = o(1) \quad \mathrm{a.s.} 
\]
Because $\pi_Q (i)$ is proportional to $\pi_{\wh{Q}} (i)/q_i$, on the preceding event, we have that
\[
	\pi_Q = \frac{\pi_{\wh{Q}}(i)\, (1+\delta_i)}{\sum_{k \in [n]} \pi_{\wh{Q}}(k)\, (1+\delta_k)} = (1+o(1))\, \pi_{\wh{Q}} (i).
\]
Since $\|\pi_{\wh{Q}} - u\|_{\mathrm{TV}} = o(1)$ a.s., we find that
\[
	\| \pi_Q - u \|_{\mathrm{TV}} \leq \| \pi_{\wh{Q}} - u \|_{\mathrm{TV}} + o(1) = o(1) \quad \mathrm{a.s.}
\]
\end{proof}


\noindent 
\textbf{Acknowledgements.} We thank two anonymous reviewers for their valuable feedback. The work in this paper was supported by the National Science Foundation (NSF) under grant no.\ DMS-1928930 while the authors were in residence at the Simons Laufer Mathematical Sciences Institute in Berkeley, California, during the Spring 2025 semester. J.C.\ and D.R.\ further acknowledge support through the NSF award CCF-2106687 and the US Army Research Office Multidisciplinary University Research Initiative award W911NF-19-1-0233.



\newcommand{\etalchar}[1]{$^{#1}$}


\end{document}